\documentclass[reprint,aps,prx,twocolumn,longbibliography,floatfix]{revtex4-2}

\usepackage[T1]{fontenc}
\usepackage[utf8]{inputenc}
\usepackage{microtype}
\usepackage{amsmath,amssymb,amsfonts}
\usepackage{booktabs}
\usepackage{graphicx}
\usepackage{hyperref}
\usepackage{amsthm}
\usepackage{quantikz}

\DeclareMathOperator{\Ran}{Ran}

\newcommand{\C}{\mathbb{C}}
\newcommand{\norm}[1]{\left\lVert #1 \right\rVert}

\newtheorem{assumption}{Assumption}
\newtheorem{definition}{Definition}
\newtheorem{theorem}{Theorem}
\newtheorem{proposition}{Proposition}
\newtheorem{corollary}{Corollary}
\newtheorem{remark}{Remark}

\hypersetup{
    colorlinks=true,
    linkcolor=blue,
    urlcolor=blue,
    citecolor=blue,
    hypertexnames=false
}

\begin{document}

\title{McLachlan-projected reduced dynamics for ill-posed Schr\"odingerized backward diffusion}

\author{Jeongbin Jo}
\email{jeongbin033@yonsei.ac.kr}
\affiliation{Department of Physics, Yonsei University, Seoul 03722, Republic of Korea}

\date{\today}

\begin{abstract}
Backward diffusion is a prototype ill-posed evolution: high spatial frequencies grow exponentially in time, so mesh-based time marching without explicit regularization is quickly overwhelmed by noise. Schr\"odingerization embeds the semidiscrete generator into Hermitian dynamics on an extended space; we ask whether McLachlan projection onto a fixed low-dimensional frame supplies a structured regularizer whose error budget can be read from a projection defect that separates full lifted propagation from the reduced trajectory. We prove uniqueness of the reduced flow, Gram-norm conservation, a lifted--reduced gap bound in terms of that defect, and perturbation estimates that highlight overlap-matrix conditioning when matrix elements are estimated statistically. We also spell out a fair classical baseline---spectral low-pass or Tikhonov filtering on the same semidiscrete model, with bandwidth or ridge strength matched to the information content of the chosen frame---so numerical contrasts isolate the Schr\"odingerized reduced pipeline rather than an unregularized Crank--Nicolson march that mainly showcases blow-up. A calibrated one-dimensional benchmark pairs a spectrally truncated reference with snapshot-built subspace evolution and finite-shot Qiskit Aer estimation, illustrating how lift, projection, and sampling layers contribute differently to the overall error.
\end{abstract}

\maketitle

\section{Introduction}
\label{sec:intro}

The Schr\"odingerization program converts linear ordinary and partial differential equations with nonunitary dynamics into Schr\"odinger-type systems in one higher dimension through a warped-phase transformation, thereby permitting Hamiltonian simulation of problems that are not unitary in their original variables~\cite{Jin2024quantum,Jin2023applications,Jin2022technical}.  Subsequent work has extended that framework to physical boundary and interface conditions, inhomogeneous linear systems, explicit circuit constructions, and ill-posed problems such as backward heat equations~\cite{Jin2024boundary,Hu2024quantumcircuits,Jin2025inhomogeneous,Jin2025stable}.  In particular, the recent inhomogeneous and ill-posed analyses emphasize that recovery from the extended state is itself a nontrivial algorithmic layer, not a cosmetic post-processing step~\cite{Jin2025inhomogeneous,Jin2025stable}.

This paper concerns the next layer once a lift is available: whether every time query must propagate the entire lifted amplitude vector, or whether a fixed-dimensional reduced pencil suffices. Krylov pipelines already estimate small matrices for spectral summaries~\cite{Cortes2022quantumkrylov,Zhang2024measurementefficient,Shen2023realtimekrylov}; here the same overlaps and Hamiltonian matrix elements furnish the generator of \emph{time-evolving} coefficients in the span of chosen lifted snapshots~\cite{Sirovich1987snapshotPOD}.

Backward reconstruction of diffusion fields appears in inverse heat conduction~\cite{Beck1985inverse} and deconvolution pipelines; quantum Schr\"odingerization formulations for unstable generators are now discussed in the open literature~\cite{Jin2024quantum,Jin2025stable}, but comparative numerics rarely separate \emph{lifting}, \emph{subspace truncation}, and \emph{overlap estimation}. This work supplies that separation for a fixed McLachlan frame.

The contribution stays narrow:
\begin{enumerate}
    \item We fix the motivating PDE as semidiscrete \emph{backward} diffusion (ill-posed forward-in-time evolution) and explain why any mesh integrator must pair with an explicit regularizer; we specify classical spectral and Tikhonov filters that serve as fair baselines once their cutoff is matched to the frame dimension~$m$.
    \item We isolate post-lift reduced dynamics behind Assumption~\ref{ass:lift}, derive the McLachlan equation $iS\dot c=H_Kc$, Gram conservation, and the lifted--reduced gap through $\delta_V$ (Theorem~\ref{thm:fixed-frame-error}).
    \item We concatenate lift, projection, and sampled-pencil layers in Corollary~\ref{cor:layered-error} and quantify generator perturbations through $\kappa(S)$ (Proposition~\ref{prop:sampling}).
\end{enumerate}
We do \emph{not} claim that Schr\"odingerization removes ill-posedness in the recovered physical variable, nor that projection supersedes optimal Tikhonov design. The point is methodological: projection error is an interpretable knob that plays the same conceptual role as a chosen low-pass bandwidth, enabling like-for-like classical comparisons instead of unstable Crank--Nicolson-only reference curves.

\section{Problem setting and lifted realization}
\label{sec:problem}

\subsection{Semidiscrete PDE}

Let $X$ be the state space of a linear PDE,
\begin{equation}
    \partial_t u(t;\mu)=\mathcal{A}(\mu)u(t;\mu),
    \qquad
    u(0;\mu)=u_0(\mu),
\label{eq:continuous-pde}
\end{equation}
where $\mu$ is a parameter vector.  After a spatial discretization with resolution parameter $h$, write
\begin{equation}
    \dot u_h(t;\mu)=A_h(\mu)u_h(t;\mu),
    \qquad
    u_h(0;\mu)=u_{0,h}(\mu),
\label{eq:semidiscrete}
\end{equation}
with $u_h\in\C^{N_h}$.  The present analysis is agnostic to the choice of finite differences, finite elements, or spectral discretization; all PDE-specific information enters through $A_h$ and the spatial discretization error.  Let $\mathcal I_h:\C^{N_h}\to X$ be a bounded reconstruction operator and define
\begin{equation}
    \varepsilon_{\mathrm{space}}(h,t;\mu)
    =
    \norm{u(t;\mu)-\mathcal I_hu_h(t;\mu)}_X.
\label{eq:space-error}
\end{equation}
Throughout, $\norm{\cdot}_2$ denotes the Euclidean vector norm and its induced spectral matrix norm on finite-dimensional spaces.

\subsection{Backward diffusion and ill-posedness}

For constant $\nu>0$, consider \emph{backward} parabolic evolution on $(0,L)$,
\begin{equation}
    \partial_t u(t,x)
    =
    -\nu\,\partial_{xx}u(t,x),
    \qquad
    x\in(0,L),
\label{eq:backward-heat}
\end{equation}
supplemented by linear boundary conditions. Let $D_h\preceq0$ denote the usual symmetric interior Laplacian (second difference scaled by grid spacing); forward diffusion is $\dot u_h=\nu D_h u_h$ with damping high modes. The semidiscrete \emph{backward} model tracked here is sign-reversed,
\begin{equation}
    \dot u_h(t)
    =
    A_h u_h(t),
    \qquad
    A_h := -\nu D_h\succeq0,
\label{eq:backward-semidiscrete}
\end{equation}
whose eigenvalues are nonnegative. For sinusoidal modes of wavenumber $k$, coefficients grow like $e^{+\nu k^2 t}$ (mesh-dependent): arbitrarily fine spatial structure is amplified, hence the continuum limit is ill-posed and discrete simulations require explicit \emph{regularization}~\cite{Engl1996regularization}.

\paragraph*{Fair classical filters on the same mesh.}
Two stabilizers on $\mathbb{C}^{N_h}$ are (i)~evolution that applies backward growth only to modes with $|k|\le k_\star$ after diagonalizing $A_h$ (or multiplying Fourier amplitudes by a sharp window in periodic geometry), and (ii)~\emph{Tikhonov}-stabilized steps~\cite{Tikhonov1977solutions} that penalize energy in $\operatorname{Ran}(D_h)$. Choosing $k_\star$ or the ridge parameter from a target smoothness scale ties the classical curve to an \emph{information budget} comparable to selecting an $m$-dimensional snapshot/Krylov span: both discard high-frequency energy and trade bias for stability. Unfiltered Crank--Nicolson or direct matrix-exponential evaluation for $A_h$ alone are \emph{not} fair performance references; they illustrate blow-up of noise and should appear only as sanity checks beside matched filters.

Inhomogeneous terms can be appended,
\begin{equation}
    \dot u_h
    =
    A_h u_h
    +
    b_h(t),
\label{eq:inhomogeneous-semidiscrete}
\end{equation}
after which Schr\"odingerized lifts require the extended-state recipes for affine dynamics~\cite{Jin2025inhomogeneous}.

\subsection{Schr\"odingerized realization}

Rather than rederive a particular warped-phase construction, we isolate the exact interface needed by the reduced layer.

\begin{assumption}[Lift--recovery realization]
\label{ass:lift}
For each $(h,\mu)$ and auxiliary resolution parameter $q$, there exist a lifted dimension $M_{h,q}$, a Hermitian matrix
$H_{h,q}(\mu)\in\C^{M_{h,q}\times M_{h,q}}$, a bounded injection
$J_{h,q}:\C^{N_h}\to\C^{M_{h,q}}$, and a bounded recovery map
$R_{h,q}:\C^{M_{h,q}}\to\C^{N_h}$ such that, for $0\le t\le T$,
\begin{equation}
    \begin{aligned}
    &\norm{
    u_h(t;\mu)
    -
    R_{h,q}\,e^{-itH_{h,q}(\mu)}\,J_{h,q}u_{0,h}(\mu)
    }_2
    \\
    &\qquad\le
    \varepsilon_{\mathrm{lift}}(h,q,t;\mu).
    \end{aligned}
\label{eq:lift-realization}
\end{equation}
\end{assumption}

\noindent
The Schr\"odingerization constructions already established in the literature provide realizations of this form for homogeneous linear systems, for boundary/interface problems after augmentation, and for inhomogeneous or ill-posed systems with suitable recovery choices in the extended variable~\cite{Jin2024quantum,Jin2023applications,Jin2024boundary,Jin2025inhomogeneous,Jin2025stable}.  The important point for the present work is that the reduced algorithm does not alter the lift; it acts on the Hermitian propagator supplied by Assumption~\ref{ass:lift}.

For readability, suppress $(h,q,\mu)$ and write $M\equiv M_{h,q}$ for the lifted-space dimension.  Define
\begin{equation}
    \Psi_0 = J u_{0,h},
    \qquad
    \Psi(t)=e^{-itH}\Psi_0.
\label{eq:lifted-state}
\end{equation}
The target quantity is the recovered state $R\Psi(t)$, which approximates $u_h(t)$ up to $\varepsilon_{\rm lift}$.

\section{McLachlan-projected lifted dynamics}
\label{sec:projected}

We now fix a frame of $m$ lifted snapshots and seek the best approximation to the full lifted dynamics within their span.  The guiding principle is McLachlan's variational principle~\cite{McLachlan1964variational}: among all trial states $\Psi_m(t)=Vc(t)$ in the span of a fixed frame $V$, we choose the coefficient trajectory $c(t)$ that minimizes the norm of the residual by which the Schr\"odinger equation $i\dot\Psi=H\Psi$ fails to be satisfied.  This yields a reduced Hermitian evolution whose generator is determined by the overlap and Hamiltonian matrices projected onto the frame.

\subsection{Projected frame and pencil}

Let $V=[v_1,\ldots,v_m]\in\C^{M\times m}$ have full column rank (here $M$ is the lifted-space dimension defined above), and define
\begin{equation}
    S = V^\dagger V,
    \qquad
    H_K = V^\dagger HV,
    \qquad
    P_V = V S^{-1}V^\dagger.
\label{eq:pencil-projector}
\end{equation}
Then $S\succ0$, and $P_V$ is the orthogonal projector onto $\Ran(V)$.  A reduced trajectory is written as
\begin{equation}
    \Psi_m(t)=Vc(t).
\label{eq:reduced-ansatz}
\end{equation}

\begin{definition}[McLachlan residual~\cite{McLachlan1964variational}]
For a fixed frame $V$, define
\begin{equation}
    \mathcal{J}(\dot c;c)
    =
    \norm{V\dot c+iHVc}_2^2.
\label{eq:mclachlan-functional}
\end{equation}
\end{definition}

\begin{theorem}[Projected coefficient dynamics]
\label{thm:projected-dynamics}
For fixed full-rank $V$ and Hermitian $H$, the unique minimizer of
Eq.~(\ref{eq:mclachlan-functional}) with respect to $\dot c^\ast$ satisfies
\begin{equation}
    iS\dot c = H_K c.
\label{eq:projected-ode}
\end{equation}
Equivalently,
\begin{equation}
    i\dot\Psi_m = P_VH\Psi_m,
\label{eq:projected-state}
\end{equation}
and the residual obeys the Galerkin orthogonality condition
\begin{equation}
    V^\dagger(\dot\Psi_m+iH\Psi_m)=0.
\label{eq:galerkin-orthogonality}
\end{equation}
\end{theorem}

\noindent
\emph{Remark.}  Since $S\succ0$, $\mathcal{J}$ is strictly convex in $\dot c$, so the stationary point $\partial\mathcal{J}/\partial\dot c^\ast=0$ is the unique global minimum.  The equivalence of Eqs.~\eqref{eq:projected-ode} and~\eqref{eq:projected-state} follows by multiplying the former by $VS^{-1}$, and the Galerkin condition~\eqref{eq:galerkin-orthogonality} by left-multiplying the state residual by $V^\dagger$.

\begin{corollary}[Gram-norm conservation]
\label{cor:gram-conservation}
For fixed $V$, the projected Hermitian flow preserves the Gram norm:
\begin{equation}
    \frac{d}{dt}\bigl(c^\dagger Sc\bigr)=0.
\end{equation}
If $V$ is orthonormal, then $S=I$ and Eq.~(\ref{eq:projected-ode}) reduces to ordinary Schr\"odinger evolution on $\C^m$ generated by the Hermitian matrix $H_K$.
\end{corollary}

\noindent
\emph{Remark.}  Differentiating and using $iS\dot c=H_K c$ with $H_K^\dagger=H_K$ gives $\dot c^\dagger Sc + c^\dagger S\dot c = ic^\dagger H_K c - ic^\dagger H_K c = 0$.

\begin{remark}[Relation to classical linear projection]
If one applies the same construction directly to a nonunitary semidiscrete system
$\dot u_h=A_hu_h$, then minimizing $\norm{V\dot c-A_hVc}^2_2$ yields
$S\dot c=A_Kc$ with $A_K=V^\dagger A_hV$.  Equation~(\ref{eq:projected-ode}) is the Hermitian post-lift special case obtained by taking $A_h=-iH$.
\end{remark}

\section{Error structure}
\label{sec:error}

We now quantify how far the reduced trajectory deviates from the full lifted state, and how statistical estimation of the pencil $(S,H_K)$ propagates into the recovered physical variable.  The key quantity is the \emph{projection defect}, which measures the component of the Hamiltonian action that lies outside the chosen frame.

\subsection{Projection defect}

\begin{definition}[Projection defect]
\label{def:projection-defect}
For the reduced trajectory $\Psi_m(t)$, define
\begin{equation}
    \eta_V(t)
    =
    \norm{(I-P_V)H\Psi_m(t)}_2.
\label{eq:projection-defect}
\end{equation}
\end{definition}

\begin{definition}[Cumulative projection error]
\label{def:cumulative-projection}
For a fixed frame $V$, define
\begin{equation}
    \delta_V(t)
    =
    \norm{(I-P_V)\Psi_0}_2
    +
    \int_0^t \eta_V(s)\,ds.
\label{eq:cumulative-projection-error}
\end{equation}
\end{definition}

\noindent
Since $\Psi_m(t)$ solves the smooth ODE~\eqref{eq:projected-ode}, the map $t\mapsto\eta_V(t)$ is continuous and hence integrable on $[0,T]$.

\begin{remark}[Quantifying $\eta_V$]
For orthonormal Krylov matrices one recovers $\eta_{V_m}=|h_{m+1,m}|\,|e_m^\top c|$ from the last Arnoldi coefficient. Practical implementations often build $V$ from snapshots $\exp(-iH\tau_\ell)\Psi_0$ instead; then $\eta_V$ is evaluated from the explicit residual $(I-P_V)H\Psi_m$.
\end{remark}

\begin{theorem}[Recovered-state error for a fixed frame]
\label{thm:fixed-frame-error}
Let $\Psi(t)$ solve Eq.~(\ref{eq:lifted-state}), and let $\Psi_m(t)$ solve
Eq.~(\ref{eq:projected-state}) with initial condition
$\Psi_m(0)=P_V\Psi_0$.  Then, for $0\le t\le T$,
\begin{equation}
    \norm{\Psi(t)-\Psi_m(t)}_2
    \le
    \delta_V(t).
\label{eq:lifted-error-bound}
\end{equation}
Consequently, under Assumption~\ref{ass:lift},
\begin{equation}
    \norm{u_h(t)-R\Psi_m(t)}_2
    \le
    \varepsilon_{\mathrm{lift}}(h,q,t;\mu)
    +
    \norm{R}_2\,\delta_V(t).
\label{eq:recovered-error-bound}
\end{equation}
\end{theorem}

\noindent
\emph{Remark.}  Let $e(t)=\Psi(t)-\Psi_m(t)$.  Then $i\dot e = He + (I-P_V)H\Psi_m$ with $e(0)=(I-P_V)\Psi_0$.  By the variation-of-constants formula and unitarity of $e^{-itH}$,
\begin{equation}
    \norm{e(t)}_2
    \le
    \norm{e(0)}_2
    +
    \int_0^t \norm{(I-P_V)H\Psi_m(s)}_2\,ds
    = \delta_V(t).
\end{equation}
Equation~\eqref{eq:recovered-error-bound} follows from $u_h-R\Psi_m=[u_h-R\Psi]+R[\Psi-\Psi_m]$ and Assumption~\ref{ass:lift}.

\subsection{Sampled pencil perturbations}

In practice, the entries of $S$ and $H_K$ may be estimated statistically (e.g., by Monte Carlo or Hadamard tests on a quantum device), introducing sampling noise.  We now bound the resulting perturbation in the reduced generator and the propagated coefficients.

Suppose hardware or Monte Carlo estimation returns
\begin{equation}
    \widehat S = S+\Delta S,
    \qquad
    \widehat H_K = H_K+\Delta H.
\end{equation}
Define the exact and sampled coefficient generators
\begin{equation}
    G = S^{-1}H_K,
    \qquad
    \widehat G = \widehat S^{-1}\widehat H_K.
\end{equation}

\begin{proposition}[Sampled-pencil perturbation]
\label{prop:sampling}
If $\norm{S^{-1}\Delta S}_2<1$, then
\begin{equation}
    \norm{\widehat G-G}_2
    \le
    \frac{\norm{S^{-1}}_2}
    {1-\norm{S^{-1}\Delta S}_2}
    \left(
    \norm{\Delta H}_2
    +
    \norm{\Delta S}_2\,\norm{G}_2
    \right).
\label{eq:generator-perturbation}
\end{equation}
Moreover, if $c$ and $\widehat c$ evolve from the same initial data under
$i\dot c=Gc$ and $i\dot{\widehat c}=\widehat G\widehat c$, then
\begin{equation}
    \norm{\widehat c(t)-c(t)}_2
    \le
    t\sqrt{\kappa(S)\kappa(\widehat S)}\,
    \norm{\widehat G-G}_2\,
    \norm{c(0)}_2.
\label{eq:coefficient-perturbation}
\end{equation}
\end{proposition}

\noindent
\emph{Remark.}  Since $\widehat G-G=\widehat S^{-1}(\Delta H-\Delta SG)$ and the Neumann series gives $\norm{\widehat S^{-1}}_2\le\norm{S^{-1}}_2/(1-\norm{S^{-1}\Delta S}_2)$, Eq.~\eqref{eq:generator-perturbation} follows.  For Eq.~\eqref{eq:coefficient-perturbation}, writing $S=L^2$ with $L=S^{1/2}$, the similarity transform $L^{-1}GL^{-1}=L^{-1}H_K L^{-1}$ is Hermitian, so $\norm{e^{-iGt}}_2\le\kappa(S)^{1/2}$.  Duhamel's principle then yields the stated bound with $\kappa(S)^{1/2}\kappa(\widehat S)^{1/2}$ from the two propagators.

\begin{remark}[Why orthonormalization matters]
Equation~(\ref{eq:generator-perturbation}) makes the conditioning issue explicit:
small overlap perturbations are harmless only when $\norm{S^{-1}}_2$ is controlled.
Arnoldi or Lanczos orthonormalization fixes $S=I$ in exact arithmetic, removing this amplification channel before shot noise is considered.
\end{remark}

\subsection{Total error decomposition}

Let $\widehat\Psi_m(t)=V\widehat c(t)$ be the sampled reduced lifted state and define
\begin{equation}
    \varepsilon_{\mathrm{sample}}(t)
    =
    \norm{V}_2\,\norm{\widehat c(t)-c(t)}_2.
\end{equation}
For a fixed frame, set
\(
\varepsilon_{\mathrm{proj}}(t)=\delta_V(t)
\).

\begin{corollary}[Layered physical error bound]
\label{cor:layered-error}
Define the recovered sampled approximation
\(
\widetilde u_m(t)=\mathcal I_hR\widehat\Psi_m(t)
\).
Then
\begin{multline}
    \norm{u(t)-\widetilde u_m(t)}_X
    \le
    \varepsilon_{\mathrm{space}}(h,t)
    +
    \norm{\mathcal I_h}_2
    \bigg(
    \varepsilon_{\mathrm{lift}}(h,q,t)
    \\
    +
    \norm{R}_2
    \big(
    \varepsilon_{\mathrm{proj}}(t)
    +
    \varepsilon_{\mathrm{sample}}(t)
    \big)
    \bigg).
\label{eq:total-error}
\end{multline}
\end{corollary}

\begin{remark}[Ill-posed reconstruction and projection bias]
For $A_h\succeq0$ as in Eq.~(\ref{eq:backward-semidiscrete}), idealized noise-free data still leave the inverse problem ill-conditioned: small high-frequency components dominate the recovered signal. Theorem~\ref{thm:fixed-frame-error} separates the McLachlan price $\delta_V$ from lift and sampling errors. Interpreted on backward diffusion, increasing the effective rank of $V$ lowers $\delta_V$ but widens the spectral support that is propagated, mirroring the variance--bias trade-off in classical TSVD or Fourier cutoffs~\cite{Hansen2010DIP}. Conversely, a small $m$ suppresses noise amplification at the expense of structured bias. The ``natural filter'' viewpoint is therefore qualitative: $P_V HP_V$ does not implement an optimal Wiener filter unless $V$ is chosen for that statistical objective.
\end{remark}

\noindent
Equation~(\ref{eq:total-error}) is the main conceptual payoff of introducing a reduced layer after the lift: it separates spatial discretization, lift/recovery, projection, and sampled-pencil effects that would otherwise be collapsed into one empirical discrepancy.
Table~\ref{tab:layer-accounting} collects the dominant symbols beside order-of-magnitude cost cues---intended as a design checklist rather than rigid complexity theorems.

\begin{table*}[t]
  \centering
  \footnotesize
  \setlength{\tabcolsep}{4pt}
  \caption{\label{tab:layer-accounting}%
  Layer-wise error bookkeeping tied to~\eqref{eq:total-error} alongside representative classical and quantum-facing resource scalings ($N_h$: physical dof; $M$: lifted dimension after Schr\"odingerization; $m$: reduced frame dimension; $\mathrm{nnz}$: sparse matrix nonzeros).
  Sampling row quotes the generator perturbation from Proposition~\ref{prop:sampling}; depths depend on synthesis and hardware.}
  \begin{tabular}{@{}p{0.15\textwidth} p{0.28\textwidth} p{0.26\textwidth} p{0.26\textwidth}@{}}
    \toprule
    Stage & Error / conditioning & Classical cost cues & Quantum / estimator cues \\
    \midrule
    Spatial mesh
    & $\varepsilon_{\mathrm{space}}$ [\eqref{eq:space-error}]
    & Mesh storage $O(N_h)$; one semidiscrete step costs $O(\mathrm{nnz}(A_h))$ for sparse generators
    & Optional block encodings tied to $\mathcal{A}_h$, not enumerated here \\
    \addlinespace[2pt]
    Lift--recovery
    & $\varepsilon_{\mathrm{lift}}$ [\eqref{eq:lift-realization}]
    & Dense lifted propagation naive $O(M^3)$ per matrix exponential; Hamiltonian solves use structure of $H$
    & Hamiltonian simulation depth for $e^{-itH}$, ancillas from auxiliary encoding index~$q$ \\
    \addlinespace[2pt]
    McLachlan projection
    & $\delta_V$, $\eta_V$ [\eqref{eq:cumulative-projection-error}], \eqref{eq:recovered-error-bound} after recovery
    & Reduced flow $iS\dot c=H_K c$ advances an $m\times m$ Hermitian pencil; dense linear algebra typically $O(m^3)$ per direct diagonalization timestep (cheaper splitting if applicable)
    & Once $(\widehat S,\widehat H_K)$ are assembled, coefficient evolution mimics classical $m$-dim Schr\"odinger solve \\
    \addlinespace[2pt]
    Exact pencil assembly
    & anchors $\varepsilon_{\mathrm{proj}}=\delta_V$ before sampling
    & $m$ matrix--vectors $H v_j$, then $m(m+1)/2$ nontrivial overlaps to fill $S$ and Hermitian $H_K$
    & Pauli transition machinery / Hadamard tests analogous to~\cite{Cortes2022quantumkrylov,Zhang2024measurementefficient}; count scales with decomposition of $H$ \\
    \addlinespace[2pt]
    Sampled overlaps
    & $\varepsilon_{\mathrm{sample}}$; $\|\widehat G-G\|_2$ bounded by~\eqref{eq:generator-perturbation}; $\kappa(S)=\|S\|_2\,\|S^{-1}\|_2$
    & Same linear algebra structures if estimating classically via Monte Carlo
    & Statistical error often scales as $\varepsilon^{-2}$ shots per amplitude;~\eqref{eq:generator-perturbation} shows $\|S^{-1}\|_2$ amplification \\
    \addlinespace[2pt]
    Fair classical filters
    & bias from cutoff $k_\star$ / Tikhonov parameter (paired to frame bandwidth on backward diffusion)
    & Periodic or circulant surrogates admit $O(N_h\log N_h)$ spectral filtering; ridge solves reuse sparse $D_h$
    & Not intrinsic to projector path; retained as matched classical budget reference \\
    \bottomrule
  \end{tabular}
\end{table*}

\section{Estimator cost and resource accounting}
\label{sec:hardware}

Table~\ref{tab:layer-accounting} already itemizes pencil assembly versus sampling amplification; here we emphasize the estimator interface alone.

Classical assembly of $(S,H_K)$ costs $m$ matrix--vector applies of $H$ to frame vectors and $O(m^2)$ overlaps. If
\(
H=\sum_{\alpha} h_\alpha P_\alpha
\)
and $|v_j\rangle=U_j|\phi\rangle$, the overlaps and Pauli transition amplitudes familiar from quantum Krylov recipes estimate the requisite matrix elements~\cite{Cortes2022quantumkrylov,Zhang2024measurementefficient}. Compiled circuit depth depends strongly on synthesis choice (naive Hamiltonian exponentials versus low-order Lie--Trotter); representative depth tables accompany the archived numerical workflows noted under Data Availability.

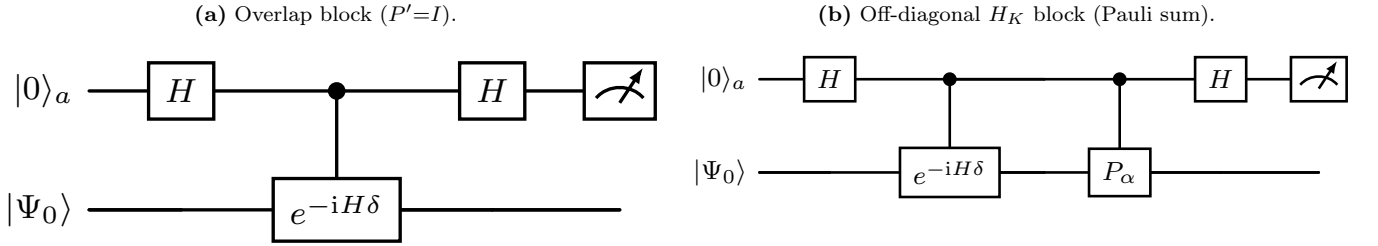
\begin{figure*}[t]
  \centering
  \footnotesize
  \setlength{\belowcaptionskip}{6pt}
  \begin{minipage}[t]{0.49\textwidth}
    \centering
    \textbf{(a)} Overlap block ($P'{=}I$).\\[0.35em]
    \resizebox{\linewidth}{!}{%
    \begin{quantikz}[column sep=0.5cm]
    \lstick{$|0\rangle_a$} & \gate{H} & \ctrl{1} & \gate{H} & \meter{} \\
    \lstick{$|\Psi_0\rangle$} & \qw & \gate{\ensuremath{e^{-\mathrm{i} H\delta}}} & \qw & \qw \\
    \end{quantikz}%
    }
  \end{minipage}\hfill
  \begin{minipage}[t]{0.49\textwidth}
    \centering
    \textbf{(b)} Off-diagonal $H_K$ block (Pauli sum).\\[0.35em]
    \resizebox{\linewidth}{!}{%
    \begin{quantikz}[column sep=0.48cm]
    \lstick{$|0\rangle_a$} & \gate{H} & \ctrl{1} & \qw & \ctrl{1} & \gate{H} & \meter{} \\
    \lstick{$|\Psi_0\rangle$}
      & \qw
      & \gate{\ensuremath{e^{-\mathrm{i} H\delta}}}
      & \qw
      & \gate{\ensuremath{P_{\alpha}}}
      & \qw
      & \qw \\
    \end{quantikz}%
    }
  \end{minipage}
  \caption{Ancillary $|0\rangle_a$ Hadamard-test schematics for estimating reduced pencil entries after Schr\"odingerization. \textbf{(a)} Controlled evolution $e^{-\mathrm{i}H\delta}$ with $P'{=}I$, supplying real and imaginary parts of $\bra{\Psi_0}e^{-\mathrm{i}H\delta}\ket{\Psi_0}$ for Gram (overlap) reconstruction. \textbf{(b)} Controlled $e^{-\mathrm{i}H\delta}$ followed by controlled Pauli $P_\alpha$ from $H=\sum_\alpha h_\alpha P_\alpha$, feeding off-diagonal projected Hamiltonian entries. Imaginary quadrature inserts $S^\dagger$ before the final Hadamard on the ancilla (Qiskit's \texttt{sdg}). Diagonal matrix elements $(H_K)_{\ell\ell}$ follow from separate snapshot expectation-value estimates after state preparation.}
  \label{fig:hadamard_schematic}
\end{figure*}

\section{Numerical Experiments}
\label{sec:numerics}

To empirically validate McLachlan projection as structured regularization and to probe estimator effects from Proposition~\ref{prop:sampling}, we simulate 1D semidiscrete backward diffusion with $\nu=0.05$, $N_x{=}7$ interior unknowns ($L{=}1$ in code), auxiliary grid $N_p{=}16$ with lifted dimension $8 N_p{=}128$, final time $T{=}0.3$, and finite-shot estimation through Qiskit's \texttt{aer\_simulator} backend at $10^4$ shots ($546$ Hadamard-test sampler circuits in the schematic of Fig.~\ref{fig:hadamard_schematic} plus $4$ estimator calls for diagonal entries).

We define the spectral reference trajectory by classical spectral truncation retaining the $k_\star{=}4$ slowest-growing modes on the homogeneous interior Laplacian (matched to the frame dimension),
and instantiate an $m{=}4$ snapshot frame composed of lifted states \(\exp(-\mathrm{i}H\tau_\ell)\Psi_0\) with \(\tau_\ell\in\{0.0,\,0.1,\,0.2,\,0.3\}\).

\subsection{Regularization by Projection}

As measured against the spectral truncation at $T{=}0.3$ ($\Delta t{=}0.03$, $11$ time samples), an unregularized Crank--Nicolson stencil reaches a $\approx\!41.99\%$ discrete relative $L^2$ endpoint error (concurrent with the plotted maximum along $t$). Evolving the augmented lifted state exactly limits the mismatch to $\varepsilon_{\mathrm{lift}} \approx 19.37\%$. An $m{=}4$ McLachlan frame with classical pencils (\emph{Projected exact matrices}) gives $\approx 20.02\%$ at $T$, whereas Aer-sampled pencils (\emph{Projected quantum matrices}) terminate at $\approx 19.82\%$ despite a larger transient excursion with peak relative error $\approx 22.03\%$ before $T$. All reduced trajectories stay in the tens-of-percent regime while Crank--Nicolson departs past $40\%$, underscoring that mesh marching without a bandwidth-matched filter is not a fair reference for ill-posed growth.

\begin{figure*}[t]
    \centering
    \begin{minipage}[t]{0.48\linewidth}
      \centering
      \textbf{(a)}\\[0.25em]
      \includegraphics[width=\linewidth]{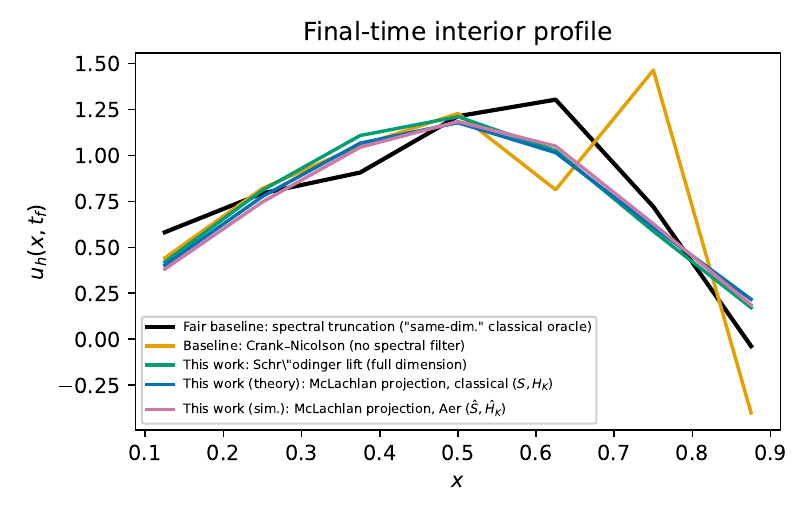}
    \end{minipage}\hfill
    \begin{minipage}[t]{0.48\linewidth}
      \centering
      \textbf{(b)}\\[0.25em]
      \includegraphics[width=\linewidth]{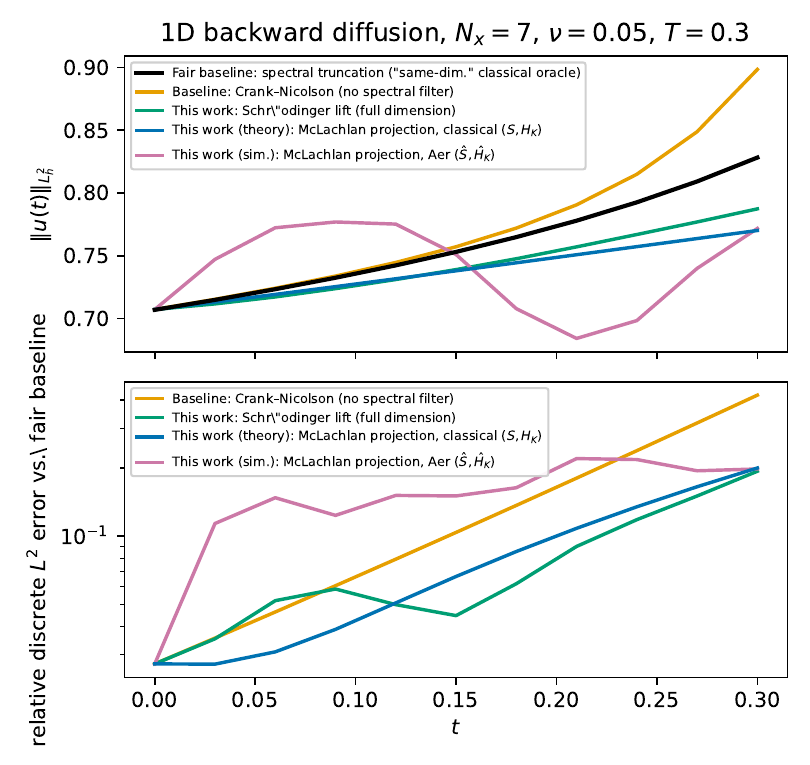}
    \end{minipage}
    \caption{\textbf{(a)} Recovered interior states $u_h$ at $T{=}0.3$. \textbf{(b)} Discrete relative $L^2$ error versus the spectral low-pass baseline ($\Delta t{=}0.03$). Crank--Nicolson attains $\approx\!41.99\%$ at the final sample (also the plotted maximum over $t$). Lifted exact propagation ends at $\approx\!19.37\%$, classical McLachlan pencils at $\approx\!20.02\%$, and Aer-sampled pencils at $\approx\!19.82\%$ with a transient peak $\approx\!22.03\%$, qualitatively consistent with~\eqref{eq:generator-perturbation}.}
    \label{fig:numerics_combined}
\end{figure*}

\subsection{Masking Hamiltonian Shot Noise}

The practical viability of quantum subspace algorithms hinges on their tolerance to sampling noise $\varepsilon_{\mathrm{sample}}$. For our $m{=}4$ frame, the snapshot basis $V$ is far from orthogonal, yielding a highly ill-conditioned overlap matrix $S$ with $\kappa(S) \approx 1.727\times 10^4$.

When $\widehat{H}_K$ is assembled from the same finite-shot estimator, shot-noise perturbations disproportionately perturb singular directions of $S$ with small overlaps. Consequently, naive matrix-wise comparisons in ambient Frobenius norm report $\widehat{H}_K$ discrepant from $H_K$ despite mild overlap errors---here $\|\widehat{H}_K-H_K\|_F/\|H_K\|_F \approx 7.08$.

Critically, this large matrix-space residual does \emph{not} manifest as exponential blow-up in the recovered observable: the reconstructed quantum trajectory in Fig.~\ref{fig:numerics_combined} remains in the same tens-of-percent regime as lifted propagation, with endpoint error $\approx 19.82\%$ and peak error $\approx 22.03\%$ before $T$ (\(\kappa(\widehat S)\approx 4.18\times 10^2\) for the sampled overlap; $\|\widehat S-S\|_F/\|S\|_F\approx 4.57\times 10^{-3}\)).

For the \emph{same} stored pair $(S,\widehat S)$, the operator norm $\|S^{-1}(\widehat S-S)\|_2\approx 5.1\times 10^{1}$ exceeds unity, so the sufficient Neumann hypothesis $\|S^{-1}\Delta S\|_2<1$ in Proposition~\ref{prop:sampling} is \emph{not} met for the ideal Gram matrix even though the Frobenius relative error is small---exactly the $\|S^{-1}\|_2$ amplification highlighted in~\eqref{eq:generator-perturbation}.

The time-stepping driver documented with the code applies the same eigenvalue floor at $10^{-2}$ to both overlaps before solving for the reduced generator; writing $S_\star,\widehat S_\star$ for those spectrally floored matrices, the reproduced numerics give $\|S_\star^{-1}(\widehat S_\star-S_\star)\|_2\approx 5.8\times 10^{-1}<1$, so the perturbation step used in the proof applies to the pencil actually integrated even though it fails for the raw $S$.

Exact reduced evolution preserves $c^\dagger Sc$ (Corollary~\ref{cor:gram-conservation}); under sampled pencils the dynamics no longer satisfy that conservation law verbatim, yet propagated coefficients \(\widehat c(t)\) still avoid the catastrophic blow-up hinted at by \(\|\widehat H_K-H_K\|_F\) alone when read out of dynamical context.

Recovery and projection bias still separate all reduced curves from the spectral reference and from unfiltered Crank--Nicolson.
Rather than asserting optimality or ``proof'' of masking, these numbers illustrate the amplification predicted by~\eqref{eq:generator-perturbation}: measured \(\Delta H\) remains large relative to \(\|H_K\|_2\) (Frobenius ratio $\approx 7.08$), \(\Delta S\) comparatively small (Frobenius ratio $\approx 4.57\times 10^{-3}$), yet the recovered endpoint errors remain in the $10\%$--$40\%$ range rather than mirroring raw Hamiltonian perturbation blow-by-blow.

\section{Discussion}
\label{sec:discussion}

The McLachlan projector supplies a deterministic reduced Hermitian evolution $iS\dot c=H_K c$ atop any validated Schr\"odingerized lift (Table~\ref{tab:layer-accounting} summarizes where each error amplitude and coarse resource budget enters). Embedding backward diffusion reframes projection loss $\delta_V$ as one explicit knob in an ill-posed inverse: enlarging span$(V)$ increases propagated bandwidth (variance grows) whereas shrinking $m$ enforces smoothing (bias grows). Classical spectral or Tikhonov filters enact the same trade-off on mesh data; juxtaposing those curves with reconstructed $RVc(t)$ separates lift artifacts, estimator noise, and the structured bias inherited from frame design.

Filtering also highlights why unconditional claims about Schr\"odingerized unitarity curing physical blow-up mislead: Eq.~(\ref{eq:recovered-error-bound}) tracks the observable after recovery maps. Fair comparisons pit matched regularizers against one another; unfiltered time marching---quantum or classical---merely diagnoses instability.

Outstanding work aligns Morozov-style discrepancy selection for classical cutoff parameters~\cite{Morozov1984methods,Hansen2010DIP} with data-driven picks of snapshots defining $V$, documents resource scalings alongside linear-combination Hamiltonian simulation~\cite{An2023linear}, and reinstates reproducible circuits once backward-heat drivers land in software.

\begin{acknowledgments}
The author acknowledges the use of IBM Quantum services for the Qiskit Aer simulations that informed the numerical experiments. The views expressed are those of the author, and do not reflect the official policy or position of IBM or the IBM Quantum team.
\end{acknowledgments}

\bibliography{bibliography}

\end{document}